\documentclass[12pt,leqno,draft]{article}

\newtheorem{theorem}{Theorem}
\newtheorem{lemma}[theorem]{Lemma}
\newtheorem{proposition}[theorem]{Proposition}
\newtheorem{definition}[theorem]{Definition}
\newtheorem{corollary}[theorem]{Corollary}

\newcommand{\begintheorem}{\addtocounter{equation}{1}\begin{theorem}}
\newcommand{\beginlemma}{\addtocounter{equation}{1}\begin{lemma}}
\newcommand{\beginproposition}{\addtocounter{equation}{1}\begin{proposition}}
\newcommand{\begindefinition}{\addtocounter{equation}{1}\begin{definition}}
\newcommand{\begincorollary}{\addtocounter{equation}{1}\begin{corollary}}

\begin{document}

\title{Notes on normed algebras, 4}

\author{Stephen William Semmes	\\
	Rice University		\\
	Houston, Texas}

\date{}

\maketitle

	Let $\mathcal{A}$ be a finite-dimensional algebra over the
complex numbers with nonzero identity element $e$.  If $x \in
\mathcal{A}$, then the resolvent set associated to $x$ is the set
$\rho(x)$ of complex numbers $\lambda$ such that $\lambda \, e - x$ is
invertible, and the spectrum of $x$ is the set $\sigma(x)$ of complex
numbers $\lambda$ such that $\lambda \, e - x$ is not invertible.  For
instance, if $V$ is a finite-dimensional vector space over the complex
numbers of positive dimension, $\mathcal{L}(V)$ is the algebra of
linear operators on $V$, and $T$ is a linear operator on $V$, then a
complex number $\lambda$ lies in the spectrum of $T$ if and only if
$\lambda \, I - T$ has a nontrivial kernel.  This is equivalent to
saying that there is a nonzero vector $v \in V$ such that $T(v) =
\lambda \, v$, which is to say that $v$ is a nonzero eigenvector for
$T$ with eigenvalue $\lambda$.

	Let $p(z)$ be a polynomial on the complex numbers, which can
be written explicitly as
\begin{equation}
	p(z) = c_m \, z^m + c_{m-1} \, z^{m-1} + \cdots + c_0
\end{equation}
for some complex numbers $c_0, \ldots, c_m$.  If $\mathcal{A}$ is a
finite-dimensional algebra over the complex numbers with identity
element $e$ and $x$ is an element of $\mathcal{A}$, then we can define
$p(x)$ in the obvious manner, by
\begin{equation}
	p(x) = c_m \, x^m + c_{m-1} \, x^{m-1} + \cdots + c_0 \, e.
\end{equation}
Here $x^j$ denotes the product of $x$ with itself with a total of $j$
copies of $x$, and we can think of this as being equal to the identity
element $e$ when $j = 0$.

	The action of polynomials on $\mathcal{A}$ behaves nicely with
respect to sums, products, and compositions, in the sense that if
$p_1$, $p_2$ are two polynomials and $x \in \mathcal{A}$, then $(p_1 +
p_2)(x) = p_1(x) + p_2(x)$, $(p_1 \, p_2)(x) = p_1(x) \, p_2(x)$, and
$(p_1 \circ p_2)(x) = p_1(p_2(x))$.  A nonzero polynomial $p(z)$ on
the complex numbers can be factored, which is to say that it can be
expressed as the product of a nonzero complex number and finitely many
factors of the form $(z - \zeta)$, where the $\zeta$'s are complex
numbers which are the zeros of $p$, counted with appropriate
multiplicity.  If $x \in \mathcal{A}$, then we can write $p(x)$ as the
product of the same complex number times $(x - \zeta \, e)$ using the
same complex numbers $\zeta$.  Hence $p(x)$ is invertible in
$\mathcal{A}$ if and only if $(x - a \, e)$ is invertible for each
complex number $\zeta$ which is a zero of $p$.

	If $p(z)$ is a polynomial, $\mu$ is a complex number, and $x$
is an element of $\mathcal{A}$, it follows that $\mu$ is in the
resolvent set of $p(x)$ if and only if $\lambda$ is in the resolvent
set of $x$ for each complex number $\lambda$ such that $p(\lambda) =
\mu$.  Therefore $\mu$ lies in the spectrum of $p(x)$ if and only if
there is a complex number $\lambda$ in the spectrum of $x$ such that
$p(\lambda) = \mu$.  In short, $p$ maps the spectrum of $x$ onto
the spectrum of $p(x)$.

	By an involution on $\mathcal{A}$ we mean an operation $x
\mapsto x^*$ with the following properties.  First, the involution is
conjugate-linear, which is to say that $(x + y)^* = x^* + y^*$ for all
$x, y \in \mathcal{A}$ and $(\lambda \, x)^* = \overline{\lambda} \,
x^*$ for all complex numbers $\lambda$ and all $x \in \mathcal{A}$,
where $\overline{\lambda}$ denotes the complex conjugate of the
complex number $\lambda$.  Second, $(x \, y)^* = y^* \, x^*$ for all
$x, y \in \mathcal{A}$.  Third, $(x^*)^* = x$ for all $x \in
\mathcal{A}$.

	For instance, one can think of the complex numbers themselves
as a one-dimensional commutative algebra, and complex conjugation
defines an involution on the complex numbers.  More generally, if $X$
is a nonempty finite set, then we can define an involution on the
algebra of complex-valued functions on $X$ by taking the complex
conjugate of such a function.  Now suppose that $V$ is a
finite-dimensional complex vector space of positive dimension equipped
with a nondegenerate Hermitian form $\langle \cdot, \cdot \rangle$.
This means that for each $v, w \in V$ we get a complex number $\langle
v, w \rangle$ such that $v \mapsto \langle v, w \rangle$ is a linear
functional on $V$ for each $w \in V$, $\langle w, v \rangle$ is equal
to the complex conjugate of $\langle v, w \rangle$ for all $v, w \in
V$, and for each nonzero $v \in V$ there is a nonzero $w \in V$ such
that $\langle v, w \rangle \ne 0$.  In this situation for each linear
operator $T$ on $V$ there is a unique linear operator $T^*$ on $V$,
called the adjoint of $T$, such that $\langle T(v), w \rangle$
is equal to $\langle v, T^*(w) \rangle$ for all $v, w \in V$,
and $T \mapsto T^*$ defines an involution on the algebra of linear
operators on $V$.

	If $\mathcal{A}$ is an algebra with involution $*$ and
identity element $e$, then the involution maps $\mathcal{A}$ onto
itself, and maps nonzero elements of $\mathcal{A}$ to nonzero elements
of $\mathcal{A}$, since $(x^*)^* = x$ for all $x \in \mathcal{A}$.  In
particular, $e^*$ is a an identity element of $\mathcal{A}$, and
therefore $e^* = e$.  If $x$ is an element of $\mathcal{A}$, and if
$x$ is invertible in $\mathcal{A}$, then $x^*$ is also invertible,
with $(x^*)^{-1} = (x^{-1})^*$.  As a result, a complex number
$\lambda$ lies in the resolvent set of $x$ if and only if the complex
conjugate $\overline{\lambda}$ of $\lambda$ lies in the resolvent set
of $x^*$.  Similarly, $\lambda$ lies in the spectrum of $x$ if and
only if $\overline{\lambda}$ lies in the spectrum of $x^*$.

	Let us call an operation $x \mapsto x^*$ on a
finite-dimensional algebra $\mathcal{A}$ with nonzero identity element
$e$ a binvolution if it satisfies the same properties as an involution
except that the multiplicative property is replaced by $(x \, y)^* =
x^* \, y^*$ for all $x, y \in \mathcal{A}$.  This is the same as an
involution when $\mathcal{A}$ is a commutative algebra.  One can check
that the properties of an involution mentioned in the preceding
paragraph work just as well for a binvolution.  If one looks more
precisely at left or right inverses in a noncommutative algebra, then
some differences appear in that an involution switches between left
and right inverses while a binvolution keeps them the same.

	Here is a basic family of examples of binvolutions.  Fix a
positive integer $n$, and consider the vector space ${\bf C}^n$ of
$n$-tuples of complex numbers.  Linear transformations on ${\bf C}^n$
can be identified with matrices in the usual manner, with respect to
the standard basis.  We can define a binvolution acting on the algebra
of linear transformations on ${\bf C}^n$ by simply taking the complex
conjugates of the entries of the corresponding matrices.
By contrast, for the usual adjoint of such a linear transformation,
corresponding to the standard inner product on ${\bf C}^n$,
one takes the complex conjugate of the transposes of the matrices.

	Suppose that $\mathcal{A}$ is a finite-dimensional algebra
over the complex numbers with nonzero identity element $e$ and an
involution or binvolution $x \mapsto x^*$.  An element $x$ of
$\mathcal{A}$ is said to be self-adjoint if $x^* = x$.  The
self-adjoint elements of $\mathcal{A}$ form a real vector subspace of
$\mathcal{A}$, and every element of $\mathcal{A}$ can be decomposed
uniquely as the sum of a self-adjoint element of $\mathcal{A}$ and $i$
times a self-adjoint element of $\mathcal{A}$.  The product of two
commuting self-adjoint elements of $\mathcal{A}$ is also self-adjoint.
In the case of an involution, $x \, x^*$ is always self-adjoint for
any $x \in \mathcal{A}$.

	Fix a self-adjoint element $x$ of $\mathcal{A}$.  Thus $x^j$
is also self-adjoint for all positive integers $j$.  The subalgebra of
$\mathcal{A}$ generated by $x$, $e$ is the same as the linear span of
$e$ and $x^j$ as $j$ runs through the positive integers.  We can
restrict the involution or binvolution $*$ to this subalgebra, which
is to say that if $y$ is in the subalgebra, then so is $y^*$.  This
subalgebra is commutative, and so it does not matter whether we
started with an involution or a binvolution on $\mathcal{A}$, in the
sense that the restriction to the subalgebra will be both.

	In the case of a normed algebra $(\mathcal{A}, \|\cdot \|)$,
there is a very natural compatibility condition to consider between
the norm and an involution or binvolution $*$, which is that
$\|x^*\| = \|x\|$ for all $x \in \mathcal{A}$.

	Let us look at a very interesting class of examples, namely
the convolution algebra on a finite semigroup.  Specifically, let $A$
be a finite semigroup with identity element $\theta$, and let
$\mathcal{A}$ be the space of complex-valued functions on $A$, as a
vector space.  If $f_1$, $f_2$ are two functions on $A$, then their
convolution $f_1 * f_2$ is the function on $A$ defined by saying that
$(f_1 * f_2)(a)$ is equal to the sum of $f_1(b) \, f_2(c)$ over all
$b, c \in A$ such that $a = b, c$.  Clearly the convolution of two
functions is linear in the functions, and one can check that it is
also associative.

	For each element $a$ of $A$, let $\delta_a$ denote the
function on $A$ which is equal to $1$ at $a$ and to $0$ at all other
elements of $A$.  Because $\theta$ is the identity element for $A$,
one can check that $\delta_\theta$ acts as an identity element for
convolution, which is to say that if $f$ is a function on $A$, then
$\delta_\theta * f$ and $f * \delta_\theta$ are equal to $f$.  Thus
$\mathcal{A}$ becomes an algebra with nonzero identity element using
convolution.  There is a simple way to define a binvolution on
$\mathcal{A}$, which is to take the complex conjugate of a function on
$A$.  If $A$ is a group, so that every element of $A$ has an inverse,
then we can define an involution on $\mathcal{A}$ by associating to a
function $f(a)$ on $A$ the function $\overline{f(a^{-1})}$.


\begin{thebibliography}{4}


\bibitem {1} W.~Arveson, {\it A Short Course on Spectral Theory},
Springer-Verlag, 2002.

\bibitem {2} R.~Beals, {\it Topics in Operator Theory},
University of Chicago Press, 1971.

\bibitem {3} S.~Berberian, {\it Introduction to Hilbert Space},
Oxford University Press, 1961.

\bibitem {4} E.~Lorch, {\it Spectral Theory}, Oxford University
Press, 1962.



\end{thebibliography}
\end{document}